\documentclass[10pt]{amsart}
\usepackage{amsmath, amsthm, amssymb}
\usepackage{mathscinet}
\usepackage{graphicx}
\usepackage{mathbbol}
\usepackage{txfonts}
\usepackage[usenames]{color}
\usepackage{enumerate}
\usepackage{hyperref}

\newtheorem{thm}{Theorem}[section]

\newtheorem{defn}[thm]{Definition}
\newtheorem{rmrk}[thm]{Remark}   

\newcommand{\be}{\begin{equation}}
\newcommand{\ee}{\end{equation}}
\newcommand{\bee}{\begin{equation*}}
\newcommand{\eee}{\end{equation*}}

\newcommand{\N}{\mathbb{N}}

\newcommand{\R}{\mathbb{R}}

\newcommand{\disjointunion}{\sqcup}

\newcommand{\injrad}{\operatorname{injrad}}

\newcommand{\diam}{\operatorname{Diam}}

\newcommand{\vol}{\operatorname{Vol}}
\newcommand{\Ricci}{\rm{Ricci}}
\newcommand{\secc}{\rm{sec}}

\newcommand{\GHto}{\stackrel { \textrm{GH}}{\longrightarrow} }

\newcommand{\dil}{\operatorname{dil}}

\begin{document}

\title{A survey on the Convergence of Manifolds with Boundary}

\author{Raquel Perales}
\thanks{Doctoral student at Stony Brook. Partially supported by 
NSF DMS 10060059.}
\address{SUNY at Stony Brook}
\email{praquel@math.sunysb.edu}

\maketitle

\begin{abstract}
This survey reviews precompactness theorems for classes of Riemannian manifolds with boundary.
We begin with the works of Kodani, Anderson-Katsuda-Kurylev-Lassas-Taylor and Wong. We then present new results of Knox 
and the author with Sormani.
\end{abstract}


\section{Introduction}

Given a sequence of Riemannian manifolds, one can say that the sequence converges if the manifolds are resembling more and more some metric space. This survey reviews theorems that are able to tell when some sequences of manifolds have subsequences
that converge, using a variety of notions of convergence presented in Section ~\ref{TypesofC}. Section \ref{sect-smoothB} is devoted to those theorems that
require smooth boundary conditions and Section~ \ref{sect-PerSor} states theorems
for open manifolds for which the boundary is defined to be $\partial M : = \bar{M}\setminus M$. 

Section ~\ref{TypesofC} reviews the definitions of Lipschitz, $C^{1,\alpha}$ and Gromov-Hausdorff distance. In addition, Cheeger-Gromov, Gromov and Anderson precompactness theorems for Riemannian manifolds without boundary are stated (Theorems ~\ref{thm-GromovLip},~\ref{GH-compactness-1},~\ref{BGcomparison},~\ref{thm-Anderson1}). Although there are interesting theorems proven both for manifolds with and without boundary for the intrinsic flat distance \cite{SorWen1,SorWen2} and for ultralimits \cite{DW-ultralimit,Gromov-asympt}, these are not discussed here.

Section ~\ref{sect-smoothB} surveys results concerning Riemannian manifolds that have smooth boundaries. 

In 1990, Kodani \cite{Kodani-1990} proves a precompactness theorem for Riemannian manifolds with boundary using the Lipschitz topology. He  assumes that the manifolds have uniformly bounded sectional curvature, 
nonnegative second fundamental forms uniformly bounded from above, and a uniform lower bound for the volume of the manifolds. See Subsection ~\ref{Kodani}. In 2004, Anderson-Katsuda-Kurylev-Lassas-Taylor \cite{Anderson-2004}  proves a precompacness theorem using a different approach under completely different hypotheses. They use $C^{1,\alpha}$ convergence ($\alpha<1$) instead of Lipschitz convergence. They assume uniform bounds on the norm of the Ricci tensor of the manifolds and their boundaries, the mean curvatures, the diameter of the manifolds, and on three different radii: injectivity, interior and boundary. See Subsection ~\ref{ssect-And}.
In 2008, Wong \cite{Wong-2008} uses Alexandrov spaces to prove two precompactness theorems in the Gromov-Hausdorff topology with the same conditions as Kodani, except that the volume bound is replaced by a diameter bound. He assumes a lower bound on sectional curvatures of the manifolds, second fundamental forms and diameter of the manifolds bounded above. Furthermore, Wong also proves that the sectional curvature bound can be replaced by
a Ricci curvature bound. See Subsection ~\ref{ssect-Wong}. 
Recently, Knox with the same technique used by Anderson-Katsuda-Kurylev-Lassas-Taylor proves a precompactness theorem in the $C^\alpha$, $\alpha<1$, and $L^{1,p}$ topologies with hypothesis analogous to theirs. He assumes bounds on the secional curvatures of the manifolds and their boundaries, the mean curvatures, the diameter of the manifolds, and a uniform lower bound in the volume of the boundaries. This theorem appears in ~\ref{ssect-Knox}. 

Section~ \ref{sect-PerSor} surveys Gromov-Hausdorff precompactness theorems concerning
classes of open manifolds where no smoothness of the boundary is required. These results appear in work with Sormani with the author \cite{Persor-2013}.
Rather than proving that a sequence of manifolds converges, we study regions within the manifolds, take the Gromov-Hausdorff limits of those inner regions and then glue the limits together to create a glued limit space. We assume conditions on the Ricci curvature, volume bounded below and above, and diameter of the inner regions uniformly bounded.

I would like to thank
Monica Moreno, Fernando Galaz and Noe Barcenas, 
organizers of Taller de Vinculacion:
Matematicos Mexicanos Jovenes en el Mundo. I would like to thank my doctoral advisor, Christina Sormani, who suggested the topic of the survey and went over it, Yuri Sobolev and Ki Song who read the survey and checked the grammar.

\section{Types of convergence}\label{TypesofC}

\subsection{Lipschitz Convergence}\label{ssect-Lip}

For details about Lipschitz convergence see \cite{Gromov-metric} and \cite{BBI}. 

\begin{defn}
Let $(X, d_X)$ and $(Y, d_Y)$ be two metric spaces.
The dilation of a Lipschitz map $f:X \to Y$ is defined by 
\begin{equation*}
\dil(f)= \sup_{x \neq x'\in X} \frac{d_Y(f(x),f(x'))}{d_X(x,x')}.
\end{equation*}
\end{defn}

A function $f: X \to Y$ is called bi-Lipschitz if both $f:X \to Y$ and $f^{-1}: f(X) \to Y$ are Lipschitz maps.

\begin{defn}
The Lipschitz distance between two metric spaces 
$(X, d_X)$ and $(Y, d_Y)$ is defined by
\begin{equation*}
d_L(X,Y)= \inf_{f:X \to Y} \log(\max\{\dil(f), \dil (f^{-1})\})  
\end{equation*}
where the infimum is taken over all bi-Lipschitz homeomorphisms $f: X \to Y$.
\end{defn}

\begin{thm}[Cheeger-Gromov]\label{thm-GromovLip}
The class of connected closed m-dimensional Riemannian manifolds $M$ satisfying:
\be \label{cond-GromovLip}|\sec(M)| \leq K, \,\,\,\vol(M) \geq v
\,\,\, and \,\,\,\diam(M)\leq D
,
\ee
is precompact in the Lipschitz topology. 
\end{thm}

This theorem follows from a Gromov's precompactness theorem (Theorem 8.25 in \cite{Gromov-81a}), in which
positive uniformly bounded injectivity radii is needed, and
Cheeger's doctoral dissertation \cite{Cheeger-thesis}, which proves that the class of manifolds satisfying (\ref{cond-GromovLip}) have positive uniformly bounded injectivity radii.

\subsection{Gromov-Hausdorff Convergence}\label{ssect-Gromov}
 
More about Gromov-Hausdorff convergence can be find in \cite{Gromov-metric} and 
\cite{BBI}. Gromov's embedding theorem
appears in \cite{Gromov-81a}.
Examples and pictures about Hausdorff converging sequences can be found 
in \cite{Sor-survey}.

\begin{defn} [Hausdorff]
Let $(Z, d_Z)$ be a metric space,
the Hausdorff distance between two subsets, $A_1, A_2 \subset Z$,
is defined as
\be
d_H^Z(A_1, A_2) = \inf \Big\{ r: \, A_1 \subset T_r(A_2), \, A_2 \subset T_r(A_1)\Big\}
\ee
where the tubular neighborhood, $T_r(A)=\big\{ x\in Z: \, d_Z(x, A)<r\big\}$.
\end{defn}

Meanwhile Gromov-Hausdorff distance avoids fixing a metric space by
considering isometric embeddings of two metric spaces into a common
metric space:

\begin{defn}
A function $\varphi: (X, d_X) \to (Z, d_Z)$ between
metric spaces is an isometric embedding if 
\be
d_Z( \varphi(x_1), \varphi(x_2)) = d_X(x_1, x_2),
\ee
for all $x_1,x_2 \in X$
\end{defn}

\begin{defn} [Gromov]
Let $(X_1, d_{X_1})$ and $(X_2, d_{X_2})$ be two compact metric spaces.
The Gromov-Hausdorff distance between them 
is defined as
\be
d_{GH}\Big(\big(X_1, d_{X_1}\big),\big(X_2, d_{X_2}\big)\Big)= 
\inf\Big\{ d_Z\big(\varphi_1(X_1), \varphi_2(X_2)\big): \,\, \varphi_i: X_i\to Z \,\Big\}
\ee
where the infimum is taken over all isometric embeddings
$\varphi_i: X_i \to Z$ and all metric spaces $Z$.
\end{defn}

Gromov proved that the Gromov-Hausdorff distance is a distance on the space of isometry classes of compact metric spaces.  
In general, elements of the classes are used and the class to which they belong is never mentioned.

The most general Gromov's precompactness theorem says:
\begin{thm}\label{GH-compactness-1}[Gromov]
Let $D>0$ and $N: (0,D] \to \N$ a function. Then the
collection $\mathcal{M}^{D,N}$, of compact metric spaces $(X, d_X)$ with
 $\diam(X) \le D$ that can be covered by $N(\epsilon)$ balls of
radius $\epsilon>0$, is precompact with respect to the
Gromov-Hausdorff distance.
\end{thm}

Given $\epsilon$ and a metric space $X$, the minimum number of 
of $\epsilon$-balls needed to cover $X$ is the same as the maximum number of
pairwise disjoint $\epsilon/2$-balls in $X$. Then $N$ can be thought as 
a function that bounds the maximum number of pairwise disjoint balls of 
compact metric spaces inside a certain class. 

The converse to Theorem~\ref{GH-compactness-1} also holds.

\begin{thm}[Gromov]\label{GH-converse} 
Suppose $(X_j, d_j)$ are compact metric spaces. If
there exists $\epsilon_0>0$ such that $X_j$ contains at least $j$ disjoint
balls of radius $\epsilon_0$, then no subsequence of the $X_j$
has a Gromov-Hausdorff limit.   
\end{thm}

Thus proving precompactness theorems with respect to the
Gromov-Hausdorff distance of a certain class of compact metric
spaces is ``reduced`` to finding a function $N$ and uniform upper diameter
bound $D$. For sequences of compact Riemannian manifolds with no boundary, 
Gromov applied the Bishop-Gromov Volume Comparison Theorem
(Theorem 2.18 \cite{Persor-2013})
to control the volume of the balls and obtain the following precompactness theorem. 

\begin{thm} [Gromov]\label{BGcomparison}
Given $m\in \N$ and $D>0$, 
let $\mathcal{M}^{m,D}$ be the class of compact
$m$ dimensional
Riemannian manifolds $M$ with 
\be 
\Ricci(M) \geq 0 \,\,\,
\text{   and  }\diam(M)\le D.
\ee
Then $\mathcal{M}^{m,D}$ is precompact with respect
to the Gromov-Hausdorff distance.
\end{thm}

\subsection{$C^{1,\alpha}$ Convergence}\label{ssect-C1alpha}

For a more detailed account on $C^{k,\alpha}$ convergence, consult \cite{Petersen-text}.

\begin{defn}
Let $\{(M_i,g_i)\}$ be a sequence of $m$-dimensional Riemannian manifolds. The sequence converges
in the $C^{1,\alpha}$ topology to a $C^{1,\alpha}$ manifold $(M,g)$ if M is a $C^\infty$ manifold such that for some fixed $C^{1,\alpha}$ atlas on M compatible with its $C^\infty$ structure, $g$ is $C^{1,\alpha}$, and there are diffeomorphisms $\varphi_i: M_i \to M$, $i=1,2,3...$, 
for which $\varphi_i^*g_i \to g$ with the $C^{1,\alpha}$-norm. 
\end{defn}

\begin{rmrk}
Theorem ~\ref{thm-GromovLip} also holds for $C^{1,\alpha}$ convergence.
\end{rmrk}

$C^{1,\alpha}$ precompactness theorems for manifolds with or without boundary have been proved using the notion of $(r,N,C^{1,\alpha})$ harmonic coordinate atlas and harmonic radius. 

\begin{defn}
Let $(M,g)$ be an $m$-dimensional compact Riemannian manifold. 
$(M,g)$ has an adapted harmonic coordinate atlas $(r,N,C^{1,\alpha})$ if there exist $C>1$ and $\{B(x_k,r)) \subset M\}_{k=1}^N$ such that 
$\{B(x_k,r/2))\}_{k=1}^N$ cover $M$, $\{B(x_k,r/4))\}_{k=1}^N$ is pairwise disjoint, and for each $k$ there is an harmonic coordinate chart $u=(u_1,...,u_m):B(x_k,10r) \to \R^m$ with
\be 
 C^{-1}\delta_{ij} \leq g_{ij} \leq  C\delta_{ij}    
\ee
and
\be r^{1+\alpha}||g_{ij}(x)||_{C^{1,\alpha}}\leq C
\ee
for all $x \in B(x_k,10r)$, where $g_{i,j}=g(\nabla u_i, \nabla u_j)$ and $\nabla$ is the Levi-Civita connection.
\end{defn}

\begin{defn}
Let $(M,g)$ be a compact Riemannian manifold.
For $x \in M$, the $C^{1,\alpha}$ harmonic radius at $x$, $r_h(x)$, is the largest radius of a geodesic ball centered at $x$ for which there is a constant $C>1$ and a coordinate chart $v:B(x,r) \to \R^m$ that satisfy
\be 
 C^{-1}\delta_{ij} \leq g_{ij} \leq  C\delta_{ij}    
\ee
and
\be r^{1+\alpha}||g_{ij}(x)||_{C^{1,\alpha}}\leq C,
\ee
where $g_{i,j}=g(\nabla v_i, \nabla v_j)$.
\end{defn}

\begin{thm}[Anderson]\label{thm-Anderson1}
The class of compact, connected Riemannian $m$-manifolds $M$ satisfying
\begin{eqnarray*}
|\Ricci(M)| \leq R, & \injrad(M) \geq i\,\,\, \text{and} \,\,\,
\diam(M) \leq D, 
\end{eqnarray*}
is precompact in the $C^{1,\alpha}$ topology. 
\end{thm}

The theorem is proven by Anderson \cite{Anderson-1990} by showing first that the harmonic radii
for manifolds in this class is uniformly bounded below. Then, for given $r, N$ and $C$, he uses the fact that
the class of compact Riemannian manifolds with $(r,N,C^{1,\alpha})$ atlases is precompact in the $C^{1,\alpha'}$ topology for all $0<\alpha'<\alpha$.

\section{Precompactness Theorems for Manifolds with Smooth Boundary Conditions}\label{sect-smoothB}

\subsection{Kodani's Precompactness Theorem}\label{Kodani}

In 1990, Kodani \cite{Kodani-1990} proves a theorem in the same line as Theorem~\ref{thm-GromovLip},
except that the Riemannian manifolds that Kodani considers have boundary and
the diameter bound is replaced by bounds on the second fundamental forms of
the boundaries. See examples of neccesity of this replacement in \cite{Kodani-1990}.

\begin{thm}[Kodani]\label{thm-Kodani}\mbox{}
Given a positive integer $m$, and numbers $K,\lambda,v>0$,
the class $\mathcal{M}(m,K,\lambda,v)$ of 
connected m-dimensional Riemannian manifolds, $M$, with boundary that satisfy
\be 
|\secc(M)| \leq K, \,\,\,\ 0 \leq II \leq \lambda,\,\,\,\text{and}\,\,\,
\vol(M) \geq v,
\ee
where $II$ stands for the second fundamental form of $\partial M$, is precompact in the Lipschitz topology.
\end{thm}

The following two definitions are needed to explain the proof of Theorem \ref{thm-Kodani} and state Theorem \ref{thm-Anderson}.

\begin{defn}\label{defn-interior-inj}
Let $M$ be a Riemannian manifold with boundary and $p$ a point in its interior.
Define the interior injectivity radius of $p$, $i_{int}(p)$, to be the supremum over all $r>0$ such that if 
$\gamma:[0,t_\gamma] \to M$ is a normal geodesic with $\gamma(0)=p$, then it is minimizing from $0$ to $\min\{t_\gamma,r\}$, where $t_\gamma$ is the first time $\gamma$ intersects $\partial M$. The interior injectivity radius of M is defined as $i_{int}(M): =\inf \{ i_{int}(p)|p \in M \}$.
\end{defn}

\begin{defn}\label{defn-bound-inj}
Let $M$ be a Riemannian manifold with boundary and $p$ a point in $\partial M$.
Define the boundary injectivity radius of $p$, $i_{\partial}(p)$, to be the supremum over all $r>0$ such that there is a minimizing geodesic $\gamma:[0,r] \to M$ with $\gamma(0)=p$ normal to $\partial M$. The boundary injectivity radius of M is defined as $i_{\partial}(M): =\inf_{p\in M}{i_{int}(p)}$.
\end{defn}

Theorem \ref{thm-Kodani} is a corollary of the following:

\begin{thm}[Kodani]\mbox{}
Let $\mathcal{M}(m,K,\lambda,i)$ be the class 
of connected m-dimensional Riemannian manifolds, $M$, with boundary and
\be 
|\secc(M)| \leq K, \,\,\,\ |II| \leq \lambda,\,\,\,\,
i_{int}(M)\geq i,\,\,and\,\,i_\partial(M)\geq i,
\ee
where $II$ stands for second fundamental form of $\partial M$,
$i_{int}$ is the interior injectivity radius and $i_{\partial}$ the boundary injectivity radius. Then 
\begin{itemize}
\item
for all $\varepsilon>0$, there exists $\delta>0$ for which if $M,N \in \mathcal{M}(m,K,\lambda,i)$ and $d_{GH}(M,N)<\delta$ then
$d_{L}(M,N)<\varepsilon$.
Thus sequences in $\mathcal{M}(m,K,\lambda,i)$ that converge in Gromov-Hausdorff sense also converge in Lipschitz sense.
\item $\mathcal{M}(m,K,\lambda, v) \subset 
\mathcal{M}(m,K,\lambda,i)$
\item
$\mathcal{M}(m,K,\lambda,v)$
is precompact in the Gromov-Hausdorff topology.
\end{itemize}
\end{thm}

Proving $\mathcal{M}(m,K,\lambda, v) \subset 
\mathcal{M}(m,K,\lambda,i)$ involves finding lower bounds for 
$i_\partial$ and $i_{int}$, which is done by looking at the
conjugate radius of $M$ and the length of simple closed geodesics
in $M$, and looking at the focal radius of $\partial M$ and the length of geodesics whose endpoints are orthogonal to $\partial M$, respectively.

The fact that $\mathcal{M}(m,K,\lambda,v)$ is precompact in the Gromov-Hausdorff topology
comes from applying volume comparison theorems. First, he shows  that for all $M \in \mathcal{M}(m,K,\lambda,v)$, $\vol(M) \leq V$. Second, he shows that if $M \in \mathcal{M}(m,K,\lambda,i)$, $p \in M$ and $\varepsilon >0 $ then $\vol(B(p,\varepsilon))\leq C$ where $C>0$ is a constant that only depends on $\varepsilon$ and $K$.

Wong proves later \cite{Wong-2008} that $II$ does not have to be nonnegative and the volume condition can be replaced by a diameter condition. See Subsection ~\ref{ssect-Wong}.

\subsection{Anderson-Katsuda-Kurylev-Lassas-Taylor's Precompactness Theorem} \label{ssect-And}
This subsection reviews the precompactness theorem  
that appear in \cite{Anderson-2004} which extend Theorem 1.1 \cite{Anderson-1990} of Anderson to manifolds with boundary. 

\begin{thm}[Anderson-Katsuda-Kurylev-Lassas-Taylor]\label{thm-Anderson}
Let $\mathcal{M}(m,R,i,H_0,D)$
be the class of compact, connected Riemannian $m$-manifolds with boundary $M$ satisfying
\begin{eqnarray*}
|\Ricci(M)| \leq R, & |\Ricci(\partial M)| \leq R \\
\injrad(M) \geq i,  & i_{int}(M) \geq i, & i_b(M) \geq 2i  \\
\diam(M) \leq D, & |H|_{Lip(\partial M)} \leq H_0
\end{eqnarray*}
where 
$i_b(M)$ 
denotes the boundary injectivity radius 
of $M$, and $H$ is the mean curvature of 
$\partial M$ in $M$.
Then $\mathcal{M}(m,R,i,H_0,D)$ is precompact in the $C^{1,\alpha}$ topology. 
\end{thm}

Theorem~\ref{thm-Anderson} is proved by showing that a larger class of manifolds is precompact in the $C^{1,\alpha'}$ topology for each $0<\alpha'<1$. The second step is to show that $\mathcal{M}(m,R,i,H_0,D)$ is contained in the larger class. This part relies completely on the use of harmonic coordinates and harmonic radii for manifolds with boundary. 

\subsection{Wong's Precompactness Theorem} \label{ssect-Wong}

This theorem, which appears in \cite{Wong-2008}, is an improvement of Theorem \ref{thm-Kodani}. 
Unlike Theorem \ref{thm-Anderson}, the hypotheses do not assume
any type of injectivity radius, do not require any bound on the Ricci curvature
of the boundary, and the condition on the mean curvature vector is replaced by a condition on the second fundamental form.

\begin{thm}[Wong]\label{thm-Wong}
The class $\mathcal{M}(m,r-,\lambda^{\pm},D)$ of n-dimensional Riemannian manifolds with boundary with 
\be \Ricci(M) \geq r-, \,\,\, \lambda^- \leq II \leq \lambda^+,
\,\,\, and \,\,\,\diam(M)\leq D,
\ee
where $II$ denotes the second fundamental form of $\partial M$, is precompact in the Gromov-Hausdorff topology.
\end{thm}

The proof consists of applying Theorem~\ref{GH-compactness-1} (Gromov Compactness Theorem).
To show that the maximum number of disjoint $\varepsilon$-balls, $N(\varepsilon,M)$, for any 
$M \in \mathcal{M}(m,r^-,\lambda^{\pm},D)$ 
and $\varepsilon > 0$ is bounded above by some constant $N(\varepsilon)$.
Clearly, $N(r,M) \leq N(cr,M)$ for $c<1$. 
Wong shows that there is an isometric extension $\tilde{M}$ of $M$ that is 
an Alexandrov space. Then he proves that there are constants $c<1$ and $\tilde{D} > 0$ such that for all $M \in \mathcal{M}(m,K^-,\lambda^{\pm},D)$, $N(cr,M) \leq N(cr,\tilde{M})$
and $\diam(\tilde{M}) \leq \tilde{D}$. 
Then by volume comparison in $\tilde{M}$, $N(cr,\tilde{M}) \leq  \tilde{N}(cr)$.

\subsection{Knox's Precompactness Theorem } \label{ssect-Knox}

The following precompactness theorem appears in \cite{Knox-2012}. The approach taken
to prove it is similar to Theorem~\cite{Anderson-2004} of Anderson-Katsuda-Kurylev-Lassas-Taylor. Unlike Theorem \ref{thm-Anderson},
there are no conditions on any type of injectivity radius, but the Ricci curvature is replaced by sectional curvature and a lower bound on the volume is added. Note that 
this theorem is not an extension of an existing theorem for manifolds without boundary
because it requires a lower bound on $\vol(\partial M)$.

\begin{thm}[Knox]\label{precompactKnox}
Let $\mathcal M(m,K,H_0,D,v_\partial)$ be the class of compact connected Riemannian $m$-manifolds with connected boundary satisfying
\begin{center}
$|\sec(M)| \leq K$,  $|\sec(\partial M)| \leq K$ \\
$0 < 1/H_0 < H < H_0$\\
$\diam(M) \leq D,$ $\vol(\partial M) \geq v_\partial$,
\end{center}
where $H$ is the mean curvature.
Then $\mathcal M(m,K,H_0,D,v_\partial)$ is precompact in the $C^\alpha$ and weak $L^{1,p}$ topologies, for any $0<\alpha<1$ and any $p < \infty$.  
\end{thm}

Knox notes that if $0 < 1/H_0 < H < H_0$ is replaced by a bound on the Lipschitz norm
of $H$, then $C^2_*$ convergence can also be obtained.
The proof of Theorem~\ref{precompactKnox} follows once it is shown that $\mathcal M(m,K,H_0,D,v_\partial)$ satisfies the hypothesis of Theorem~\ref{precompactKnox2}.

\begin{thm}[Knox]\label{precompactKnox2}
If $\{(M_i,g_i)\}$ is a sequence of Riemannian manifolds with boundary such that 
\be 
r_h(g_i) \geq r_0\,\,\,\text{and}\,\,\,\diam(M_i)\leq D
\ee
where $r_h(g_i)$ is the $L^{k,p}$ harmonic radius. Then there is a subsequence of $\{(M_i,g_i)\}$ that converges in weak $L^{k,p}$ topology to a manifold with boundary whose metric is in $L^{k,p}$.
\end{thm}

The harmonic radius, $r_h(g)$, of a Riemannian manifold with boundary, $(M,g)$, depends on the harmonic radius of points in the interior of $M$ and the harmonic radius of points in $\partial M$. Knox deals with these two cases separately. First, by looking at the volume of cylinders whose base is in $\partial M$, he finds that 
there is a $c>0$ that only depends on $\mathcal M(m,K,H_0,D,v_\partial)$
such that
\be r_h(x)\geq c d_M(x, \partial M)
\ee
for all $x \in M \setminus \partial M$, where $(M,g) \in \mathcal M(m,K,H_0,D,v_\partial)$.
Then he shows that
\be r_h(x)\geq r
\ee
for all $x \in \partial M$ where $r$ is a constant that only depends on $\mathcal M(m,K,H_0,D,v_\partial)$. 
Thus, by definition of harmonic radius of a manifold with boundary, $r_h(g)$ has a lower bound that depends only on 
the class $\mathcal M(m,K,H_0,D,v_\partial)$.

\section{Precompactness Theorems for Manifolds without Boundary Conditions} \label{sect-PerSor}

This section presents a joint work between the author with Sormani appearing in \cite{Persor-2013}. We make no assumptions on the boundary. In fact, we consider open 
Riemannian manifolds, $(M^m,g)$,
endowed with the length metric, $d_M$. The boundary of $M$ is defined as
$\partial M : = \bar{M}\setminus M$
where $\bar{M}$ is the metric completion of $M$. 
Here the boundary is avoided by considering $\delta$ inner regions. Precompactness theorems are proven for these inner regions, then their Gromov-Hausdorff limit spaces are glued together
into a single metric space (Theorem~\ref{def-glue}). At the end, 
the Hausdorff dimension of this single metric space is obtained and it is shown that its Hausdorff measure has positive lower density everywhere.   
Analogous theorems for constant sectional curvature are proven
in \cite{Persor-2013}.

\subsection{$\delta$-Inner Regions and their Limits} \label{ssect-innerReg}

Let  $(M^m,g)$ be an open Riemannian manifold and $\delta>0$,
the $\delta$ inner region of $M$ is defined by 
$M^\delta:=\Big\{ x\in M: \,\, d_M(x, \partial M) > \delta\Big\}$.
There are two metrics in $M^\delta$: the restricted metric $d_M$, and the induced length metric $d_{M^\delta}$. If $M^{\delta}$ is not path connected, then the distance between points in two different path components is defined to be infinity. In general, 
$d_M(x,y)\leq d_{M^\delta}(x,y)$ for all $x,y \in M^\delta$.  

\begin{defn}
Given $m\in \N$, $\delta>0$,  $D>0$, $V>0$, and $\theta>0$, set
$
\mathcal{M}^{m, \delta, D, V}_{\theta}  
$
to be the class of m-dimensional open Riemannian manifolds, $M$, with boundary,
with 
\be 
\Ricci(M) \geq 0, \,\,\,  \vol(M)\le V, 
\text{   and  }\diam(M^\delta,{d_{M^\delta}})\le D,
\ee
that are noncollapsing at a point:
\be \label{V-epsilon}
\exists q\in M^\delta
\textrm{ such that } \vol(B_q(\delta)) \ge \theta\delta^m.
\ee
\end{defn}

\begin{thm}[P\---Sormani]\label{main-thm}
If $(M_j, g_j) \subset \mathcal{M}^{m, \delta, D, V}_{\theta}$, then
there is a subsequence $\{M_{j_k}\}$ and a compact metric
space $(Y^\delta,d)$ such that 
$
(\bar{M}^\delta_{j_k}, d_{M_{j_k}}) \GHto (Y^\delta,d_Y).
$
\end{thm}  

Note that even though $D$ bounds the diameter of the inner regions with respect to the induced length metric, the convergence is guaranteed endowing
the inner regions with the restricted metric. 

Replacing $\delta>0$ in the above theorem by a 
decreasing sequence, $\delta_i \to 0$, and adding
bounds on the diamater of $\delta_i$-inner regions the following can be proved.

\begin{thm}[P\---Sormani]\label{glued-Ricci-limits}
Given $m\in \N$, a decreasing sequence, $\delta_i \to 0$, $D_i>0$, $i=0,1,2...$, $V>0$, $\theta>0$. Suppose that  
$
\{(M_j, g_j)\} \subset \mathcal{M}^{m, \delta_0, D_0, V}_{\theta}
$
and
\be
\sup\Big\{\diam\left(M^{\delta_i}_{j}, {\,d_{M_j^{\delta_i}}}\right):\,\,j \in \N\Big\} < D_{i} \qquad
\forall i \in \N.
\ee
Then there is a subsequence $\{M_{j_k}\}$, and there are compact metric
spaces $(Y^{\delta_i},d_{Y^{\delta_i}})$ such that 
$
(\bar{M}^{\delta_i}_{j_k}, d_{M_{j_k}}) \GHto (Y^{\delta_i},d_{Y^{\delta_i}})
$ 
for all $i$.
\end{thm}

\subsection{Constructing a Glued Limit Space} \label{ssect-gluedSpace}

By constructing isometric embeddings between the limit spaces, 
$
\varphi_{\delta_{i+1}, \delta_i} : Y^{\delta_i} \to Y^{\delta_{i+1}},
$
it is possible to define a metric space into which all the limit spaces
isometrically embed. Set $\varphi_{\delta_{i+j}, \delta_i}=\varphi_{\delta_{i+j}, \delta_{i+j-1}}
\circ \cdots \circ \varphi_{\delta_{i+1}, \delta_{i}}$.
Define 
\be
Y: =Y(\{\delta_i\}, \{\varphi_{\delta_{i+1}, \delta_{i}} \} )=
Y^{\delta_0} \,\disjointunion \, \left( \,\disjointunion \,_{i=1}^\infty 
\left(Y^{\delta_{i+1}}\setminus  \varphi_{\delta_{i+1}, \delta_i}\left(Y^{\delta_i}\right)\right)\right)
\ee 
and  
\[ 
d_Y(x,y): =\left\{
\begin{array}{ll}
d_{Y^{\delta_0}}(x,y)&  \textrm{ if } x,y \in Y^{\delta_0},\\
&\\
d_{Y^{\delta_{i+1}}}(x,y)&\textrm{ if } x,y \in Y^{\delta_{i+1}}
\setminus \varphi_{\delta_{i+1}, \delta_i}\left(Y^{\delta_i}\right) ,\\
&\\
d_{Y^{\delta_{i+1}}}\left(x, \varphi_{\delta_{i+1},\delta_0}(y)\right) &
\textrm{ if } x \in Y^{\delta_{i+1}}\setminus 
\varphi_{\delta_{i+1}, \delta_{i}}\left(Y^{\delta_{i}}\right)\\
& \textrm{ for some } i \in \N \textrm{ and } y \in Y^{\delta_0}, \\
d_{Y^{\delta_{i+j+1}}}\left(x,
\varphi_{\delta_{i+j+1},\delta_{{i+1}}}(y)\right) & 
\textrm{ if } x \in Y^{\delta_{i+j+1}}\setminus 
\varphi_{\delta_{i+j+1}, \delta_{i+j}}\left(Y^{\delta_{i+j}}\right)\\
& \textrm{ and } y\in Y^{\delta_{i+1}}\setminus 
\varphi_{\delta_{i+1}, \delta_i}\left(Y^{\delta_i}\right)
\textrm{ for some } i,j \in \N 
\end{array}
\right.
\]

\begin{thm}[P\---Sormani]\label{def-glue}
Under the hypothesis of Theorem ~\ref{glued-Ricci-limits}.
There exists a metric space $(Y,d_Y)$ such that 
for all $\delta\in (0,\delta_0]$, there is a subsequence of $\{(\bar{M}^{\delta}_{j_k}, d_{M_{j_k}})\}$ that Gromov-Hausdorff converges to some compact
metric space $(Y^{\delta},d_{Y^{\delta}})$. For any such $Y^\delta$, there exists an isometric embedding 
\be
F_\delta=F_{\delta, \{\delta_i\}}:Y^{\delta} \to Y.   
\ee
If $\delta=\delta_i$ for some $i$, then 
\be \label{check-justified}
F_{\delta_i}(Y^{\delta_i}) \subset F_{\delta_{i+1}}(Y^\delta_{i+1}). 
\ee
If $\beta_j$ is any sequence decreasing to $0$, then
\be\label{eqn-onto-Y-d}
Y=\bigcup_{i=1}^{\infty} F_{\beta_j}(Y^{\beta_j}).
\ee
\end{thm}

This glued limit space may exist even when $(M_j, d_j)$
has no Gromov-Hausdorff limit.

Hausdorff measures and topologies of the Gromov-Hausdorff limit spaces of
noncollapsing sequences of manifolds have been studied by 
Cheeger, Colding, Naber, Wei and Sormani (c.f. \cite{ChCo-PartI}, 
\cite{Colding-volume}, \cite{Colding-Naber-2} and \cite{SorWei1}).
Applying some of these results, we are able to prove the following.

\begin{thm}[P\---Sormani]\label{thm-ricci-glued-lim-measure}
Suppose that $Y$ is
a glued limit constructed as in Theorem~\ref{def-glue}.
Then $Y$ has Hausdorff dimension $m$, $\mathcal{H}^m(Y)\le V_0$ and its Hausdorff measure has positive lower density everywhere.
\end{thm}

\bibliographystyle{plain}
\bibliography{2012}

\end{document}